\documentclass[11pt]{amsart}

\usepackage{amsmath,amssymb,amsthm}
\usepackage{graphicx}

\newtheorem{lemma}{Lemma}[section]
\newtheorem{theorem}[lemma]{Theorem}
\newtheorem{fact}[lemma]{Fact}

\newtheorem{definition}[lemma]{Definition}
\newtheorem{corollary}[lemma]{Corollary}
\newtheorem{proposition}[lemma]{Proposition}

\newtheorem{problem}[lemma]{Open problem}

\newcommand{\nices}{\mathcal{S}}

\newcommand{\nat}{\mathbb{N}}

\newcommand{\nicea}{\mathcal{A}}
\newcommand{\niceb}{\mathcal{B}}

\newcommand{\nicef}{\mathcal{F}}
\newcommand{\niceh}{\mathcal{H}}
\newcommand{\nicet}{\mathcal{T}}
\newcommand{\decomp}{\mathcal{R}}

\newcommand{\lk}{\textrm{lk}}

\newcommand{\sd}{\textrm{sd}}

\newcommand{\PT}{\raisebox{0.25ex}{\includegraphics{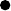}}}

\newcommand{\POne}{\raisebox{0.25ex}{\includegraphics{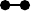}}}
\newcommand{\PTwo}{\raisebox{0.25ex}{\includegraphics{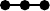}}}
\newcommand{\PThree}{\raisebox{0.25ex}{\includegraphics{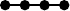}}}
\newcommand{\PFour}{\raisebox{0.25ex}{\includegraphics{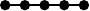}}}

\newcommand{\SThree}{\includegraphics{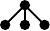}}
\newcommand{\Star}{\SThree}
\newcommand{\SpecA}{\includegraphics{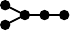}}
\newcommand{\SpecB}{\includegraphics{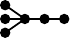}}

\begin{document}
\title[Vertex decompositions]{Vertex decompositions of two-dimensional complexes and graphs}
\author{Micha{\l} Adamaszek}
\address{Warwick Mathematics Institute and DIMAP,
      \newline University of Warwick, Coventry, CV4 7AL, UK}
\email{aszek@mimuw.edu.pl}
\thanks{Research supported by the Centre for Discrete
        Mathematics and its Applications (DIMAP), EPSRC award EP/D063191/1.}


\begin{abstract}
We investigate families of two-dimensional simplicial complexes defined in terms of vertex decompositions. They include nonevasive complexes, strongly collapsible complexes of Barmak and Miniam and analogues of 2-trees of Harary and Palmer. We investigate the complexity of recognition problems for those families and some of their combinatorial properties. Certain results follow from analogous decomposition techniques for graphs. For example, we prove that it is NP-complete to decide if a graph can be reduced to a discrete graph by a sequence of removals of vertices of degree 3.
\end{abstract}

\maketitle

\section{Introduction}
\label{section:intro}

\subsection{Overwiew.}
Simplicial complexes can be thought of as combinatorial models for topological spaces. The additional, combinatorial structure provides invariants finer than the topological ones. Many of them are defined in terms of various recursive decompositions.


Perhaps the most prominent such notion is that of \emph{collapsibility}. If a complex $K$ has a face $\sigma$ (called a \emph{free face}) which is contained in a unique maximal face $\tau$ and $\dim\tau=\dim\sigma+1$ then the removal of the pair of faces $\{\sigma, \tau\}$ is called an \emph{elementary collapse}. The complex $K$ is \emph{collapsible} if there exists a sequence of elementary collapses that reduces $K$ to the one-point complex $\PT$. This notion lies at the foundations of simple homotopy theory of Whitehead \cite{Wht}.

Most importantly an elementary collapse preserves the homotopy type of the complex (in fact it is a deformation retraction). It means that every collapsible complex is contractible. The inverse implication holds for one-dimensional complexes (a graph is collapsible, or contractible, if and only if it is a tree), but not in higher dimensions. The dunce hat of Zeeman \cite{Zee} is one of the counterexamples in dimension two. However, the existence of a collapsing sequence can always be seen as a \emph{witness} or \emph{proof} of contractibility. This matters from the complexity-theoretic point of view because it is algorithmically undecidable if a given finite simplicial complex is contractible \cite{Hak}.

There are other recursive decompositions which can produce such proofs of contractibility. In this work we concentrate on decompositions defined by removals of vertices, rather than faces. The main property of this sort is \emph{nonevasiveness}\footnote{Note that the name is slightly misleading for our purposes, because it begins with a negative prefix (non-), while in fact it represents a positive statement (the existence of a reduction).}, defined as follows. We say that a complex $K$ is nonevasive if either $K=\PT$ or there exists a vertex $v$ such that both the link $\lk_K(v)$ and the deletion $K\setminus v$ are nonevasive. Recall that the link $\lk_K(v)$ is the subcomplex of $K$ consisting of the faces $\tau\in K$ such that $v\not\in\tau$ and $\tau\cup\{v\}\in K$. Geometrically the link of a vertex $v$ is homeomorphic to the intersection of $K$ with a small sphere around $v$.

The notion of nonevasiveness arose from the topological approach of \cite{KSS} to the so-called Evasiveness Conjecture for graph properties of Karp. For more information about this see \cite[Chap.13]{Koz}. What matters for us is that the removal of a vertex $v$ as above is again a homotopy equivalence, so the existence of a nonevasive reduction proves the contractibility of $K$. 

For one-dimensional complexes (graphs) the classes of contractible, collapsible and nonevasive complexes all coincide with the class of trees. In higher dimensions the situation is a lot more interesting. Two-dimensional collapsible complexes can still be recognized by the \emph{greedy} algorithm, which collapses any of the available free faces (see e.g. \cite{HogMet}). In higher dimensions this is no longer true and the complexity of this recognition problem is an open question, although there are some NP-hardness proofs for closely related problems (\cite{EgeGon,MalFra}). For nonevasiveness the complexity status is open even in dimension two. It is likely that all these problems are computationally hard. 

\subsection{Reducibility.}

Our reference for simplicial complexes and combinatorial algebraic topology is \cite{Koz}. All complexes we consider are finite. All graphs are finite, simple, undirected and unlabeled. To avoid confusion we use the term \emph{$d$-complex} for a simplicial complex of dimension \emph{at most $d$}. A $1$-complex is simply a graph. In this paper we only deal with $2$-complexes and their $0$-, $1$- and $2$-dimensional faces are called, respectively, vertices, edges and triangles. Familiarity with the notions of geometric realization, contractibility or basic homotopy theory can aid the reader's intuition, but is not essential.

Our starting point is the observation that since nonevasive graphs are precisely trees and the link of a vertex in a $2$-complex is a graph, the definition of a nonevasive $2$-complex can be rephrased by requiring that the link $\lk_K(v)$ in each decomposition step is simply a \emph{tree}. This leads to a generalization in the form of the following notion.

\begin{definition}
Let $\nicef$ denote any family of graphs. A $2$-complex $K$ is \emph{$\nicef$-reducible} if and only if either 
\begin{itemize}
\item $K=\PT$, or
\item there is a vertex $v\in K$ such that the link $\lk_K(v)$ is a graph from $\nicef$ and $K\setminus v$ is  $\nicef$-reducible.
\end{itemize}
The family of all $\nicef$-reducible $2$-complexes is denoted $\decomp(\nicef)$.
\end{definition}

For any $2$-complex $K$ a vertex $v\in K$ such that $\lk_K(v)\in\nicef$ will be called \emph{initial} (with respect to $\nicef$) and the removal of an initial vertex will be called an \emph{elementary $\nicef$-reduction}. If there is a sequence of $\nicef$-reductions which starts with $K$ and terminates with some complex $L$ we say $K$ is $\nicef$-reducible to $L$.

We will be mostly interested in the cases when the family $\nicef$ consists of trees, because then each elementary $\nicef$-reduction removes a cone over a contractible subspace, hence it preserves the homotopy type of the complex. However, we will also mention some other types of reductions, especially in connection with graphs (see below). First let us consider some examples.


\begin{itemize}
\item If $\nicet$ denotes the class of all trees then $\decomp(\nicet)$ is, by definition, the family of all nonevasive $2$-complexes. Clearly if $\nicef_1\subset\nicef_2$ then $\decomp(\nicef_1)\subset\decomp(\nicef_2)$ so whenever $\nicef$ consists of trees then every $\nicef$-reducible $2$-complex is nonevasive.
\item The family $\decomp(\{\PT\})$ is, of course, the family of trees.
\item The family $\decomp(\{\PT, \POne\})$ is closely related to the family of \emph{2-trees} defined in \cite{HarPal}. These are the closest 2-dimensional analogues of trees. We will discuss these families in  Section \ref{section:2tree}.
\item Let $\nices=\{\PT,\POne,\PTwo,\SThree,\ldots,\SThree_n,\ldots\}$ denote the family of all \emph{stars}, that is trees of diameter at most two (the symbol $\SThree_n$ denotes a star with $n$ leaves). Then $\decomp(\nices)$ is precisely the class of \emph{strongly collapsible} $2$-complexes in the sense defined recently by Barmak and Minian in \cite{BarMin}. 
\end{itemize}

These theories are paralleled by analogous vertex-based decompositions of graphs. In a graph, viewed as a $1$-complex, the link of a vertex is a discrete space, which can be identified with its cardinality, equal to the degree of that vertex. Our previous definition can be more naturally reformulated as follows.

\begin{definition}
Let $\nicea\subset\nat$. A graph $G$ is \emph{$\nicea$-reducible} if and only if either 
\begin{itemize}
\item $G=\PT$, or
\item there is a vertex $v\in G$ whose degree is a number from $\nicea$ and such that $G\setminus v$ is $\nicea$-reducible.
\end{itemize}
\end{definition}

Of course this is a combinatorial, rather than topological notion, because vertex removals hardly ever preserve the homotopy type of the graph (this only happens when we remove a vertex of degree one). We immediately see the following reformulation. An orientation of a graph is a choice of direction for each edge. An orientation is acyclic if the resulting directed graph has no directed cycles.

\begin{fact}
A graph is $\nicea$-reducible if and only if it has an acyclic orientation such that one vertex is a sink (has out-degree zero) and the out-degree of every other vertex is in $\nicea$. In particular, if $0\in\nicea$ then a graph is $\nicea$-reducible if and only if it has an acyclic orientation such that the out-degree of every vertex is in $\nicea$. 
\end{fact}

The orientation is obtained by directing all edges out from the vertex being removed in the $\nicea$-reduction. Conversely, given an acyclic orientation we can remove any source (vertex with in-degree zero) and proceed recursively, obtaining an $\nicea$-reduction.

Certain graph properties can be expressed in these terms:
\begin{itemize}
\item Trees are the $\{1\}$-reducible graphs and forests are the $\{0,1\}$-reducible graphs.
\item The minimal $n$ such that $G$ is $\{0,1,\ldots,n\}$-reducible is usually called the \emph{degeneracy} of $G$.
\end{itemize}

\subsection{Statement of results.}

We are interested in the complexity of the membership problems for the families $\decomp(\nicef)$ and in the structure of those families. The main motivation is to understand the cases when $\nicef$ is a subclass of trees, because the associated $\nicef$-reductions preserve homotopy type and because these classes approximate the class of nonevasive $2$-complexes. For this reason we are often going to restrict attention to the most natural classes of trees $\nicef$ which are \emph{subtree-closed}, that is with each tree the class contains also all its subtrees. Some results will have equivalents for graphs.

We first use our framework to exhibit a family of trees $\nicef$ for which $\nicef$-reducibility is NP-complete and a tightly related graph-theoretic analogue.
\begin{theorem} The following problems are NP-complete for every $d\geq 3$:
\label{thm:npc}
\begin{itemize}
\item[a)] deciding if a given $2$-complex $K$ is $\{\PT,\SThree_d\}$-reducible,
\item[b)] deciding if a given graph is $\{0,d\}$-reducible.
\end{itemize}
\end{theorem}
To the author's knowledge part a) is the first example of a homotopy-type-preserving reduction for which the associated decision problem is provably hard. On the other hand, the author is aware of the fact that the families $\{\PT,\SThree_d\}$ are not as natural as one would like (e.g. not subtree-closed). Hopefully these results will be improved upon in the future.

We continue by investigating the other complexity extreme, that is the question of when $\nicef$-reducibility can be decided by the \emph{greedy} algorithm, which in each step removes any available initial vertex. This is motivated by a number of existing examples of this kind. The collapsibility of $2$-complexes can be checked greedily. The degeneracy of the graph can also be found by the greedy strategy. It follows from \cite{BarMin} that also strong collapsibility can be checked greedily and we will prove that this is essentially the only such situation for homotopy-type-preserving $\nicef$-reductions along subtree-closed families.  Let $\nices_n=\{\PT,\POne,\PTwo,\ldots,\Star_n\}$ be the family of stars with at most $n$ leaves (we allow $n=0,1,\ldots,\infty$ so that $\nices_\infty=\nices$).
\begin{theorem}
\label{prop:greedy}
Let $\nicef$ be a subtree-closed family of trees. Then membership in $\decomp(\nicef)$ can be decided by the greedy algorithm if and only if $\nicef=\nices_n$ for some $n=0,1,\ldots,\infty$.
\end{theorem}
A consequence of this theorem is that nonevasive $2$-complexes cannot be recognized greedily. As a byproduct of the proof we also obtain the next result. Recall that a family of graphs $\niceh$ is called \emph{hereditary} if it is closed under vertex removal (if $G\in\niceh$ then $G\setminus v\in\niceh$ for all $v\in G$). In particular it contains the empty graph.
\begin{proposition}
\label{prop:greedyh}
If $\niceh$ is a hereditary graph family then membership in $\decomp(\niceh)$ can be decided by the greedy algorithm.
\end{proposition}
As with graphs, this is a combinatorial, rather than topological kind of reduction. It is an immediate corollary for the family $\niceh$ consisting of all discrete graphs that the degeneracy of a graph can be computed greedily.

The argument used to prove that the greedy algorithm fails for certain families of trees can also be used to prove the following structural result.
\begin{proposition}
\label{prop:oneinitial}
Let $\nicef$ be a subtree-closed family of trees.
\begin{itemize}
\item[a)] If $\nicef=\{\PT\}$ or $\nicef=\{\PT,\POne\}$, then every complex in $\decomp(\nicef)$ has at least two initial vertices.
\item[b)] If $\nicef$ is any other subtree-closed family of trees then $\decomp(\nicef)$ contains a $2$-complex with only one initial vertex. Moreover, for any $T\in\nicef$, except possibly $T=\POne$, this $2$-complex can be chosen in such a way that the link of the unique initial vertex is isomorphic with $T$.
\end{itemize}
\end{proposition}
This result also has an expected corollary.
\begin{corollary}
If $\nicef_1\neq \nicef_2$ are two different subtree-closed families of trees then $\decomp(\nicef_1)\neq\decomp(\nicef_2)$.
\end{corollary}
We finish with brief remarks on the relation between $\decomp(\{\PT,\POne\})$ and $2$-trees and a connection between reducibility, collapsibility and nonevasiveness for barycentric subdivisions, expanding on \cite{Vel}.

\subsection*{Acknowledgement.} The author is grateful to J.A.Barmak and G.Minian for writing \cite{BarMin}, and especially for the stimulating example of \cite[Fig.7]{BarMin}.

\section{An NP-complete recognition problem}
\label{section:npc}

In this section we prove Theorem \ref{thm:npc}, starting with the graph-theoretic problem of part b).

\begin{proof}[Proof of Theorem \ref{thm:npc}, part b)]
Let us first consider the case $d=3$. We will perform a reduction from the NP-complete problem known as EXACT COVER BY 3-SETS (X3C), \cite[Prob. SP2]{Gar}.

\begin{tabular}{ll}
Instance: & A set $S$ and a family $\niceb$ of 3-element subsets of $S$.\\
Question: & Is there a subfamily $\niceb'\subset \niceb$ such that \\
 & each element of $S$ belongs to exactly one member of $\niceb'$?
\end{tabular}

Given an instance $(S,\niceb)$ of X3C we construct a graph $G$ as follows (see Fig.\ref{fig:sample}). For every element $x$ of $S$ there are three vertices called $x,x',x''$, connected by the edges $x-x'$ and $x-x''$. For every 3-set $B\in\niceb$ there is a vertex called $B$. Finally, for every pair $x,B$ such that $x\in B$ we have an edge $x-B$. Note that every vertex labeled $B$ has degree $3$, so in any orientation of $G$ which has outdegrees $0$ and $3$ such vertex will be either a source (outdegree $3$) or a sink (outdegree $0$). It also implies that any such orientation must necessarily be acyclic. Moreover, each of the edges $x-x'$, $x-x''$ must be directed outwards from $x$ (otherwise $x'$, $x''$ would have a forbidden outdegree $1$), so each $x$ must have outdegree $3$ and therefore exactly one edge directed towards one of the $B$.

\begin{figure}
\includegraphics[scale=1]{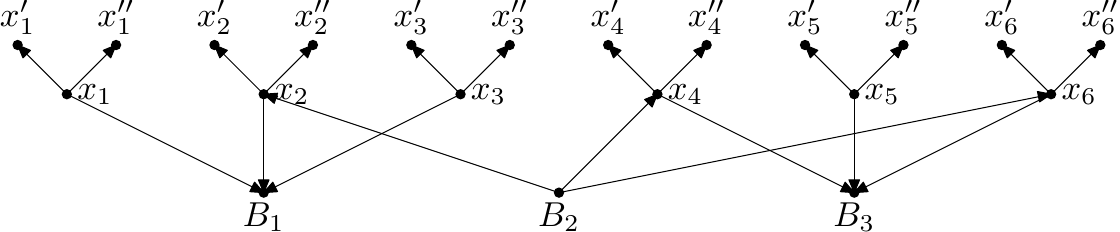}
\caption{The graph corresponding to the instance $S=\{x_1,\ldots,x_6\},B_1=\{x_1,x_2,x_3\},B_2=\{x_2,x_4,x_6\},B_3=\{x_4,x_5,x_6\}$ of X3C. The arrows indicate the acyclic orientation. The solution to the instance is $\{B_1,B_3\}$.}
\label{fig:sample}
\end{figure}

Using these remarks it is easy to see that solutions to $(S,\niceb)$ correspond to orientations of $G$ with outdegrees $0$ and $3$. More precisely, the exact cover $\niceb'$ consists of those vertices $B$ of $G$ which are sinks under the given orientation. The fact that every vertex $x$ has exactly one outgoing edge towards some 3-set $B$ corresponds to the fact that every element of $S$ belongs to exactly one 3-set from $\niceb'$.

For $d>3$ the proof is similar but uses a reduction from the NP-complete problem EXACT COVER BY $d$-SETS.
\end{proof}

Part a) of Theorem \ref{thm:npc} is deduced from the next lemma. Here $\sqcup$ denotes the disjoint union of complexes or graphs and $CK$ is the simplicial cone over $K$, i.e. the complex $CK=\{\tau,\tau\cup\{a\}: \tau\in K\}$ for a new vertex $a$ (the apex of the cone). For a subset $\nicea\subset\nat$ let $\nices_\nicea=\{\SThree_n: n\in\nicea\}$ be the family of stars whose number of leaves is in $\nicea$. Then we have the following equivalence.

\begin{proposition}
\label{dupa}
The following conditions are equivalent for a graph $G$.
\begin{itemize}
\item[a)] The $2$-complex $C(G\sqcup G)$ is $\nices_\nicea$-reducible.
\item[b)] The $2$-complex $CG$ is $\nices_\nicea$-reducible to the apex of the cone.
\item[c)] The graph $G$ is $\nicea$-reducible.
\end{itemize}
\end{proposition}
\begin{proof}
Denote by $a$ the apex of the cone in $C(G\sqcup G)$. During the $\nices_\nicea$-reduction process the link of $a$ remains disconnected (hence $a$ cannot be removed) until at least one of the two copies of $CG$ has been reduced to $a$. This proves the equivalence of a) and b). Next we prove the equivalence of b) and c). Observe that at every step of the $\nices_\nicea$-reduction of $CG$ to $a$ the link of every vertex $v$ in $G$ is the cone (with apex $a$) over the neighbourhood of  $v$ in the remaining part of $G$. This cone is $\SThree_n$ if and only if the removed vertex has degree $n$ in the remaining part of $G$. This establishes a bijection between $\nices_\nicea$-reductions of $CG$ to $a$ and $\nicea$-reductions of $G$.
\end{proof}

Part a) of Theorem \ref{thm:npc} now follows because $\nices_{\{0,d\}}=\{\PT,\SThree_d\}$.

Note that there are polynomial time algorithms recognizing the $\{0,1\}$-reducible graphs (forests) and $\{\PT,\POne\}$-reducible $2$-complexes (Theorem \ref{prop:greedy} with $\nicef=\nices_1$). The complexity of both problems in the remaining case $d=2$ is unknown.
\begin{problem}
Is there a polynomial time algorithm that decides if a graph is $\{0,2\}$-reducible? Is there a polynomial time algorithm that decides if a $2$-complex is $\{\PT,\PTwo\}$-reducible?
\end{problem}
Of course Proposition \ref{dupa} shows that the two-dimensional problem is at least as hard as the problem for graphs. The author suspects that both problems are NP-complete.

\section{Greediness and initial vertices}
\label{section:greedy}

In this section we work exclusively with the subtree-closed classes of trees. Every subtree-closed family $\nicef$ satisfies exactly one of the two conditions: either $\nicef=\nices_n$ for some $n=0,1,\ldots,\infty$ or $\{\PT,\POne,\PTwo,\PThree\}\subset\nicef$ (equivalently $\PThree\in\nicef$). We begin by analysing the greedy algorithm for $\nices_n$-reducibility and $\niceh$-reducibility for a hereditary graph class $\niceh$ (Propositions \ref{prop:greedy} and \ref{prop:greedyh}). 

Note that every star, except $\POne$, has a distinguished centre (the unique vertex of degree other than one) and we have the following simple lemma.
\begin{lemma}
\label{lemma:auto}
Let $K$ be a $2$-complex and let $v$ be a vertex of $K$. If $\lk_K(v)$ is a star other than $\POne$, whose centre is $w$, and $\lk_K(w)$ is also a star, then there is an automorphism of $K$ which swaps $v$ and $w$ and leaves all other vertices fixed.
\end{lemma}
\begin{proof}
Clearly $\lk_K(w)$ must be centered at $v$ and it follows that $v$ and $w$ belong to exactly the same simplices of $K$.
\end{proof}

\begin{proof}[Proof of correctness of the greedy algorithm for $\nices_n$-reducibility.]
Suppose that $K$ is $\nices_n$-reducible and $v_1,\ldots,v_k$ is the order in which the vertices of $K$ are removed. Suppose also that $\lk_K(v_i)\in\nices_n$ for some $i$. We will find a new order of removals which starts with $v_i$ and which is a valid $\nices_n$-reduction. This proves that the greedy algorithm can start with $v_i$ just as well.
\begin{itemize}
\item[1)]
If $\lk_K(v_i)=\POne$ or if $\lk_K(v_i)$ is a star with centre $v_j$ where $j>i$ then the new removal order is
$$v_i,v_1,\ldots,v_{i-1},v_{i+1},\ldots,v_k.$$
In the new ordering the links of $v_{i+1},\ldots,v_k$ at the time of their removal are the same as previously, while the links of $v_1,\ldots,v_{i-1}$ differ from the previous ones by at most a removal of one leaf, so they remain in $\nices_n$. It means we obtain a valid $\nices_n$-reduction.
\item[2)]
If $\lk_K(v_i)$ is a star with centre $v_j$ and $j<i$ then let $K'=K\setminus\{v_1,\ldots,v_{j-1}\}$. The link $\lk_{K'}(v_i)$ is still a star from $\nices_n$ (though perhaps smaller than $\lk_K(v_i)$) because the centre of $\lk_K(v_i)$ was never removed when passing from $K$ to $K'$. If $\lk_{K'}(v_i)=\POne$ then by 1) applied to $K'$
\begin{equation}\tag{*}\label{xx1}
v_1,\ldots,v_{j-1},v_i,v_j,v_{j+1},\ldots,v_{i-1},v_{i+1},\ldots,v_k
\end{equation}
is a valid $\nices_n$-reduction of $K$. If, on the other hand, $\lk_{K'}(v_i)$ is a star with some centre, then this centre must still be $v_j$ (it was the centre of $\lk_{K}(v_i)$ and was never removed). Moreover, $\lk_{K'}(v_j)$ is also a star (because $v_j$ is the next to be removed in an $\nices_n$-reduction). It follows by Lemma \ref{lemma:auto} that $v_i$ and $v_j$ are indistinguishable in $K'$ and
\begin{equation}\tag{**}\label{xx2}
v_1,\ldots,v_{j-1},v_i,v_{j+1}\ldots,v_{i-1},v_j,v_{i+1},\ldots,v_k
\end{equation}
is a valid $\nices_n$-reduction of $K$. In both cases (\ref{xx1}) and (\ref{xx2}) we can apply 1) to move $v_i$ to the front and complete the proof. 
\end{itemize}
\end{proof}

\begin{proof}[Proof of correctness of the greedy algorithm for $\niceh$-reducibility.]
If $\niceh$ is a hereditary graph class then the proof is much easier. Suppose that $v_1,\ldots,v_k$ is the order of vertex removals in some $\niceh$-reduction of $K$ and that some vertex $v_i$ satisfies $\lk_K(v_i)\in\niceh$. Then one can always simply bring $v_i$ forward and consider a new order
$$v_i,v_1,\ldots,v_{i-1},v_{i+1},\ldots,v_k.$$
The links of $v_1,\ldots,v_{i-1}$ during this procedure differ from the previous ones by at most a removal of one vertex. Since $\niceh$ is hereditary, those links remain in $\niceh$, so the new order defines a correct $\niceh$-reduction.
\end{proof}

\begin{figure}
\begin{tabular}{ccc}\includegraphics[scale=0.6]{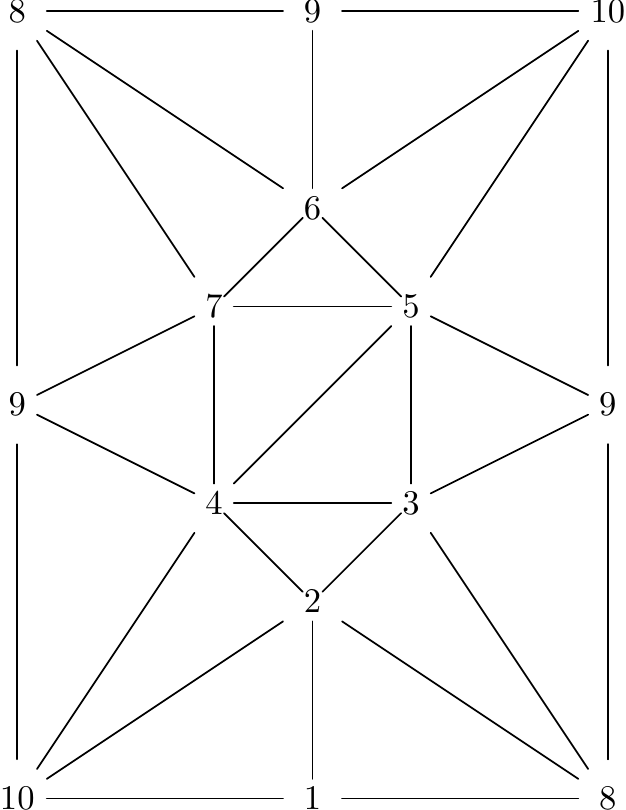} & & \includegraphics[scale=0.6]{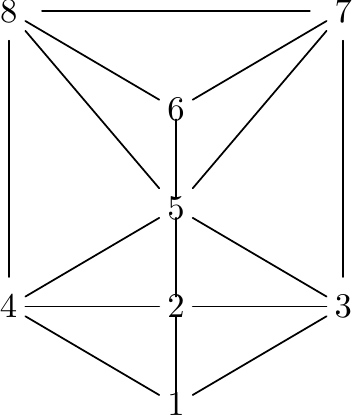}\\
a) & & b)\end{tabular}
\caption{The $2$-complexes $\mathfrak{K}$ and $\mathfrak{L}$ of Section \ref{section:greedy}. All triangles are to be filled in.}
\label{fig:ugly}
\end{figure}

\begin{figure}
\includegraphics[scale=0.6]{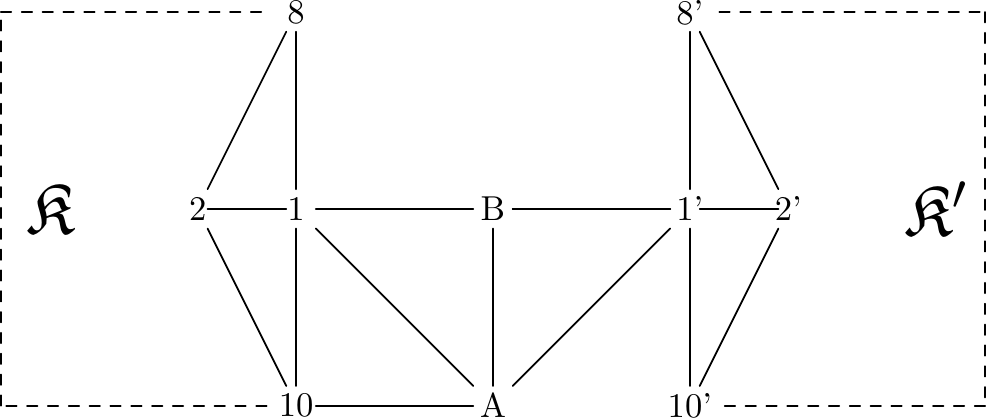}
\caption{The complex with an initial vertex $A$ which is a bad first choice.}
\label{fig:bad}
\end{figure}

\begin{proof}[Proof that the greedy algorithm for $\nicef$-reducibility fails when $\{\PT,\POne,\PTwo,\PThree\}\subset\nicef$.]
Consider the simplicial complex $\mathfrak{K}$ of Fig.\ref{fig:ugly}a. One checks by inspection that in $\mathfrak{K}$ the link of each vertex other than $1$ contains a cycle and that $\mathfrak{K}$ is $\{\PT,\POne,\PTwo,\PThree\}$-reducible to the vertex $10$ in the order of increasing vertex labels, starting with $1$.

Now consider the $2$-complex of Fig.\ref{fig:bad}, build from two copies $\mathfrak{K}$, $\mathfrak{K'}$ of $\mathfrak{K}$, two extra vertices $A,B$ and three extra triangles $\{1,10,A\}$, $\{1,A,B\}$, $\{1',A,B\}$. This complex is $\{\PT,\POne,\PTwo,\PThree\}$-reducible: first remove $B$, then reduce $\mathfrak{K}$ starting from $1$ down to $10$, remove $10$ and $A$ and reduce $\mathfrak{K}'$. Therefore the whole complex is also $\nicef$-reducible.

On the other hand the link of $A$ is $\PThree$, but if $A$ is removed first the remaining part has no initial vertex at all --- the link of each vertex is either disconnected or contains a cycle hence is not in $\nicef$. Therefore a greedy reduction algorithm will fail if it starts by choosing $A$.
\end{proof}

This proves Propositions \ref{prop:greedy} and \ref{prop:greedyh}. We continue by exploiting some further applications of $\mathfrak{K}$ and proving the second part of Proposition \ref{prop:oneinitial} (for the first part see Section \ref{section:2tree}).

\begin{proof}[Proof of Proposition \ref{prop:oneinitial}.b)]
Let $\nicef$ be any subtree-closed family other than  $\{\PT\}$ or $\{\PT,\POne\}$. In other words, suppose that $\PTwo\in\nicef$. We first exhibit a complex in $\decomp(\nicef)$ whose only initial vertex has link $\PTwo$.
\begin{itemize}
\item If $\nicef\subset\nices$ then the complex $\mathfrak{L}$ of Fig.\ref{fig:ugly}.b is $\{\PT,\POne,\PTwo\}$-reducible, hence $\nicef$-reducible. The only initial vertex is $1$, because the links of all other vertices contain either $\PThree$ or a cycle.
\item If $\PThree\in\nicef$ then the complex $\mathfrak{K}$ of Fig.\ref{fig:ugly}.a is $\{\PT,\POne,\PTwo,\PThree\}$-reducible, hence $\nicef$-reducible. The only initial vertex is $1$, because the links of all other vertices contain a cycle.
\end{itemize}
It $T=\PT$ then the complex $\mathfrak{L}$ (resp. $\mathfrak{K}$) with an extra edge $\{1,a\}$ to a new vertex $a$ is $\nicef$-reducible and $a$ is the unique initial vertex. Of course its link is $\PT$. For trees $T$ with at least three vertices we proceed by induction, which started above with $T=\PTwo$, the only 3-vertex tree. Now suppose $T$ has $n\geq 4$ vertices. Decompose $T$ as $T=T'\cup_w e$, where $e$ is an edge from a vertex called $w$ to a leaf of $T$. Since $\nicef$ is subtree-closed the tree $T'$ is in $\nicef$ and has at least three vertices, so we may assume the claim was proved for $T'$. Let $K'$ be the corresponding $\nicef$-reducible complex with a single initial vertex $v'$ satisfying $\lk_{K'}(v')=T'$. 

\begin{figure}
\begin{tabular}{ccccc}\includegraphics{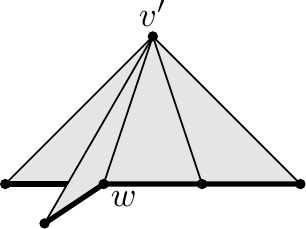} & & \includegraphics{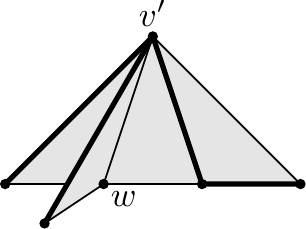} & &  \includegraphics{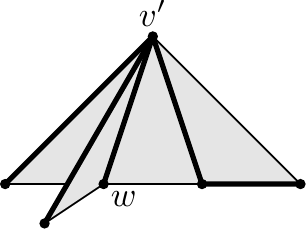}\\
a) & &  b) & & c)\end{tabular}
\caption{a) A vertex with link $T'=\protect\SpecA$. b) A modified embedding $T''$ of $\protect\SpecA$. c) Extension to an embedding of $T=\protect\SpecB$.}
\label{fig:embedding}
\end{figure}

Consider the embedding of $T'$ in $K'$ as $\lk_{K'}(v')$ and identify the vertex $w$ of $T'$ with the vertex of $\lk_{K'}(v')$ it goes to under this embedding (Fig.\ref{fig:embedding}a). Now modify the embedding so that it passes through $v'$ instead of $w$ (Fig.\ref{fig:embedding}b). Precisely, consider the tree $T''$ in the $1$-skeleton of $K'$ obtained by replacing all edges $x-w$ of $T'$ with the corresponding edges $x-v'$. This tree is isomorphic with $T'$. Now we can embed $T$ in $K'$ as $T=T''\cup_{v'}\{v'w\}$ (Fig.\ref{fig:embedding}c). The new complex $K$ is obtained by attaching a cone over this particular embedding of $T$:
$$K=K'\cup_{T}CT.$$
If $v$ denotes the apex of the cone $CT$ in $K$ then clearly $\lk_K(v)=T$. Moreover, our construction of the embedding of $T$ in $K'$ guarantees that $v'$ has degree at least $2$ as a vertex of $T$. It means that there are at least two vertices $u_1,u_2\in T'$ such that $v'u_1$, $v'u_2$ are edges of $T$. Then there is a cycle in $\lk_K(v')$ of the form
$$v-u_1-\cdots-u_2-v,$$
where $u_1-\cdots-u_2$ denotes a path connecting $u_1$ and $u_2$ in $T'$. It follows that $\lk_K(v')$ is not in $\nicef$. The links of the remaining vertices are not in $\nicef$, because
\begin{itemize}
\item if $\nicef\subset\nices$ then those links in $K'$ contained $\PThree$ or a cycle,
\item if $\PThree\in\nicef$ then those links in $K'$ contained a cycle,
\end{itemize}
and passing from $K'$ to $K$ can only enlarge the links. It means that $v$ is the only initial vertex of $K$ and it satisfies $\lk_K(v)=T$. The proof is complete.
\end{proof}

\textbf{Remark.} For certain families $\nicef$ it may be impossible to find an $\nicef$-reducible $2$-complex whose only initial vertex would have the link $\POne$. The family of all trees $\nicet$ is an example of such $\nicef$.

\begin{fact}
If $K$ is nonevasive and one of its vertices has link $\POne$ then $K$ has at least one more vertex whose link is a tree.
\end{fact}
\begin{proof}
Suppose otherwise. Let $L$ be a nonevasive complex with only one initial vertex and assume that the link of that vertex is $\POne$. By definition of nonevasiveness, the complex $L\setminus v$ is nonevasive. However, the removal of $v$ changes the links of only two other vertices, and in each of them it reduces the link by deleting from it a vertex of degree one. Such operation cannot convert a non-tree into a tree, so $L\setminus v$ has no initial vertices at all. This contradiction proves the fact.
\end{proof}


\section{$2$-trees and their relatives}
\label{section:2tree}

We now briefly discuss the class $\decomp(\{\PT,\POne\})$ and the family of $2$-trees of \cite{HarPal}. Since the interesting properties can be found in \cite{HarPal}, we will keep this description rather sketchy and omit the proofs.

A \emph{$2$-tree} is a connected $2$-complex $K$ such that every two distinct edges $e,e'$ are connected by a unique alternating walk of the form $$e=e_1-t_1-e_2-t_2-\cdots-t_{k-1}-e_k=e'$$ where $e_i$ are edges, $t_i$ are triangles, $e_i,e_{i+1}\subset t_i$ for all $i$ and all $e_i$ and $t_i$ are pairwise distinct. In particular $K$ is either pure two-dimensional or a single edge or a point (we allow the last two possibilities for convenience; in the original definition $K$ had to be pure of dimension two.)

This can also be formulated in terms of the bipartite adjacency graph $A_{1,2}(K)$ between edges and triangles of $K$ in which and edge $e$ and a triangle $t$ are adjacent if $e\subset t$. Then $K$ is a $2$-tree if it is connected and $A_{1,2}(K)$ is a tree. By looking at $A_{1,2}(K)$ one can easily verify that a $2$-tree is $\{\PT,\POne\}$-reducible, that it has at least two initial vertices with respect to $\{\PT,\POne\}$ (compare \cite[Prop.2.]{HarPal}) and that it can be reduced to any of its vertices. 

It is also easy to see that a general $\{\PT,\POne\}$-reducible complex $K$ is a simply-connected union of its maximal sub-$2$-trees such that the intersection of every two of them is either empty or a single vertex. The simply-connectedness condition can be phrased as follows. Consider the multigraph whose vertices are the maximal $2$-trees in $K$ and there is an edge for every intersection point of distinct $2$-trees. Then this multigraph must be a tree (in particular it must not have any multiple edges). Then one sees that a $\{\PT,\POne\}$-reducible complex also has at least two initial vertices and it can be reduced to any vertex.

Intuitively, the way in which $\{\PT,\POne\}$-reducible complexes are assembled from $2$-trees resembles the way in which every graph is assembled from its $2$-connected (in the sense of graph theory) components.

\section{Barycentric subdivisions}
Welker \cite{Vel} proved that the barycentric subdivision $\sd K$ of a collapsible simplicial complex $K$ is nonevasive. The following observation describes an intermediate step.

\begin{proposition}
\label{prop:subdivision}
The following conditions are equivalent for a $2$-complex $K$:
\begin{itemize}
\item[a)] $K$ is collapsible.
\item[b)] $\sd K$ is $\{\PT,\PTwo,\PFour\}$-reducible.
\item[c)] $\sd K$ is nonevasive.
\item[d)] $\sd K$ is collapsible.
\end{itemize}
\end{proposition}
Let us emphasize one corollary.
\begin{corollary}
For a $2$-complex $K$, if $\sd K$ is nonevasive, then it is in fact $\{\PT,\PTwo,\PFour\}$-reducible.
\end{corollary}
Note that not every nonevasive reduction of $\sd K$ must be of this restricted form. It would be interesting to know if for some other classes of complexes the existence of one kind of reduction implies the existence of a more restricted one. 
\begin{proof}[Proof of Proposition \ref{prop:subdivision}]
a)$\implies$b). Every elementary collapse of a free edge in $K$ can be simulated in the barycentric subdivision by the removal of a vertex with link $\PTwo$ (vertex representing an edge of $K$) followed by the removal of a vertex with link $\PFour$ (vertex representing a face of $K$). A collapse of a free vertex in $K$ is simulated by two reductions with link $\PT$ in $\sd K$.

b)$\implies$c)$\implies$d) are obvious because it is known that a nonevasive complex is collapsible.

d)$\implies$a). This follows from the known fact that for $2$-complexes collapsibility is an invariant of homeomorphism type (see e.g. \cite{HogMet}). 
\end{proof}



\begin{thebibliography}{99}

\bibitem{BarMin} J.A. Barmak, E.G. Minian, \textit{Strong homotopy types, nerves and collapses}, \texttt{arxiv/0907.2954}

\bibitem{EgeGon} \"O. E{\v g}ecio{\v g}lu, T.F. Gonzalez, \textit{A computationally intractable problem on simplicial complexes}, Comput. Geom., Theory and Applications 6 (1996) 85-98

\bibitem{Gar} M.R. Garey, D.S. Johnson, \textit{Computers and Intractability --- A Guide to the Theory of NP-Completeness}, Freeman, Oxford, UK, 1979

\bibitem{Hak} W. Haken. \textit{Connections between topological and group theoretical decision problems} In: Boone, Cannonito and Lyndon (1973), 427-441.

\bibitem{HarPal} F. Harary, E. Palmer, \textit{On acyclic simplicial complexes}, Mathematika 15 (1968), 115-122

\bibitem{HogMet} C.Hog-Angeloni, W.Metzler, \textit{Geometric aspects of two-dimensional complexes}, in \textit{Two-dimensional homotopy and combinatorial group theory}, Cambridge Univ. Press 1993, 1-50

\bibitem{KSS} J. Kahn, M. Saks, D. Sturtevant, \emph{A topological approach to evasiveness}, Combinatorica 4 (1984), 297-306

\bibitem{Koz} D. Kozlov, \textit{Combinatorial Algebraic Topology}, Algorithms and Computation in Mathematics, Vol. 21, Springer-Verlag Berlin Heidelberg 2008

\bibitem{MalFra} R. Malgouyres, A.R. Frances, \textit{Determining Whether a Simplicial 3-Complex Collapses to a $1$-Complex is NP-Complete}, Proc. DGCI 2008, LNCS 4992, 177-188, 2008

\bibitem{Vel} V. Welker, \textit{Constructions preserving evasiveness and collapsibility}, Discrete Math. 207 (1999) 243-255

\bibitem{Wht} J.H.C. Whitehead, \textit{Simple homotopy types}, Amer. J. Math. 72 (1950), 1-57

\bibitem{Zee} E.C. Zeeman, \textit{On the dunce hat}, Topology 2 (1964), 341-358
\end{thebibliography}
\end{document}